\documentclass[12pt,twoside]{article}

\usepackage{amsfonts}
\usepackage{amssymb}
\usepackage{amsmath}
\usepackage{graphicx}

\newcommand{\ignore}[1]{}

\def\@begintheorem#1#2{\par\bgroup{\sc #1\ #2. }\it\ignorespaces}
\def\@opargbegintheorem#1#2#3{\par\bgroup{\sc #1\ #2\ (#3). } \it\ignorespaces}
\def\@endtheorem{\egroup}
\newtheorem{theorem}{Theorem}[section]
\newtheorem{corollary}[theorem]{Corollary}
\newtheorem{lemma}[theorem]{Lemma}
\newtheorem{example}[theorem]{Example}
\newtheorem{proposition}[theorem]{Proposition}
\newtheorem{definition}[theorem]{Definition}
\newcommand{\bt}[1]{\begin{theorem}\label{#1}}
\newcommand{\bc}[1]{\begin{corollary}\label{#1}}
\newcommand{\bl}[1]{\begin{lemma}\label{#1}}
\newcommand{\be}[1]{\begin{example}\label{#1}}
\newcommand{\bp}[1]{\begin{proposition}\label{#1}}
\newcommand{\ba}[1]{\begin{algorithm}\rm\label{#1}}
\newcommand{\bd}[1]{\begin{definition}\rm\label{#1}}
\newcommand{\bpr}{\noindent {\em Proof. }}
\newcommand{\et}{\end{theorem}}
\newcommand{\ec}{\end{corollary}}
\newcommand{\el}{\end{lemma}}
\newcommand{\ee}{\end{example}}
\newcommand{\ep}{\end{proposition}}
\newcommand{\ed}{\end{definition}}
\newcommand{\epr}{{\ \vbox{\hrule\hbox{%
\vrule height1.3ex\hskip0.8ex\vrule}\hrule}}\\\par}
\newcommand{\mepr}{{\ \ \ \vbox{\hrule\hbox{%
\vrule height1.3ex\hskip0.8ex\vrule}\hrule}}}

\def\R{\mathbb{R}}

\def\Z{\mathbb{Z}}
\def\N{\mathbb{N}}

\def \l {\lambda}

\def \mon {{\rm mon}}

\def \size {{\rm size}}
\def \sign {{\rm sign}}

\begin{document}

\title{\bf Huge tables and multicommodity flows
are fixed-parameter tractable via unimodular integer Carath\'eodory}
\author{
Shmuel Onn
\thanks{\small Technion - Israel Institute of Technology, Haifa, Israel.
Email: onn@ie.technion.ac.il}
}

\date{}

\maketitle

\begin{abstract}
The three-way table problem is to decide if there exists an $l\times m\times n$
table satisfying given line sums, and find a table if yes.
Recently, it was shown to be fixed-parameter tractable with parameters
$l,m$. Here we extend this and show that the huge version
of the problem, where the variable side $n$ is encoded in binary, is
also fixed-parameter tractable with parameters $l,m$. We also conclude
that the huge multicommodity flow problem with a huge number of consumers is
fixed-parameter tractable. One of our tools is a theorem about
unimodular monoids which is of interest on its own right.

\end{abstract}

\section{Introduction}

The study of multiway table problems, also known as multi-index transportation
problems, goes back to the classical paper of Motzkin \cite{Mot}. It also has
applications in privacy in databases and confidential data disclosure of statistical tables,
see the survey \cite{FR} by Fienberg and Rinaldo and the references therein.
Specifically, the three-way table problem is to decide if there exists a nonnegative integer
$l\times  m\times n$ table satisfying given line sums, and find a table if there is one.
Deciding the existence of such a table is NP-complete already for $l=3$, see \cite{DO1}. Moreover, {\em every}
bounded integer program can be isomorphically represented in polynomial time for some $m$ and $n$
as some $3\times m\times n$ table problem, see \cite{DO2}. When both $l$ and $m$ are fixed, the problem
can be solved in polynomial time using Graver bases and the theory of $n$-fold integer programming
\cite{DHOW}. See the book \cite{Onn} for further background.

Recently, the problem was shown in \cite{HOR} to be fixed-parameter
tractable when $l$ and $m$ are parameters, solvable in time
$O(f(l,m)\cdot n^3\cdot\size(\mbox{line sums}))$ where
$\size(\mbox{line sums})$ is the binary-encoding length of all the given
line sums and $f(l,m)$ is a suitable computable function.
(While we do not need a bound on $f(l,m)$ for our results, it is
worthwhile to note that it is known to satisfy $f(l,m)=(lm)^{O(lm)}$,
see \cite{BO,KT,Onn} for more details on this important so-called
{\em Graver complexity} function.)
Recall that a {\em parameterized problem} with parameter $p$ and input $I$
is called {\em fixed-parameter tractable} if it admits an algorithm that runs in time
$O(f(p)\cdot\size(I)^k)$ for some computable function $f(p)$ of $p$
which is independent of $I$ and some $k$ which is independent of $p$
and $I$. In particular, if a problem is fixed-parameter tractable, then
for each fixed value $p$ of the parameter, it is polynomial-time solvable,
but fixed-parameter tractability is much stronger since the degree $k$ of
the polynomial running time is independent of the parameter value $p$.
See the book \cite{DF} by Downey and Fellows for more details on
this important branch of complexity theory.

More recently, in \cite{OS}, the huge version of the problem, where
the variable table side $n$ is a huge number encoded in {\em binary},
and the $n$ many $l\times m$ layers of the table come in
$t$ types, was also shown to be polynomial-time solvable for fixed $l$ and $m$.
Here we strengthen this and show that the huge problem is moreover fixed-parameter
tractable as well. (All layers of each given type have the same specified
row and column sums, see Section 3 for a more detailed description of the problem.)

\vskip.2cm\noindent{\bf Theorem \ref{tables}}\hskip.2cm {\em The
huge $l\times m\times n$ table problem with $t$ types, parameter
$l$, and $n$ variable and binary-encoded, is fixed-parameter
tractable in the following situations:
\begin{enumerate}
\item
when $m$ is also a parameter and $t$ is variable and unary-encoded;
\item
when $t$ is also a parameter and $m$ is variable and unary-encoded.
\end{enumerate}}

\vskip.2cm One application of this theorem is to multicommodity
flows. We show that the problem with a variable number $m$ of
suppliers and a huge variable number $n$ encoded in binary of
consumers, that come in $t$ types, is fixed-parameter tractable.
(All consumers of each given type have the same consumption in each
commodity and the same capacity from each supplier, see Section 3 for more details.)

\vskip.2cm\noindent{\bf Corollary \ref{commodity}}\hskip.2cm
{\em The huge multicommodity flow problem parameterized by the number $l$
of commodities and the number $t$ of consumer types is fixed-parameter tractable.}

\vskip.5cm
One of the tools we use is a theorem about totally unimodular monoids which we discuss next.
Let $S\subseteq\Z^d$ be a set of integer points. The {\em monoid} generated by $S$
is the set of nonnegative integer combinations of finitely many elements of $S$,
$$\mon(S)\ :=\
\left\{\sum_{k=1}^m\l_k x^k\ :\ m,\l_1,\dots,\l_m\in\Z_+\,,\ x^1,\dots,x^m\in S\right\}\ .$$
The {\em monoid decomposition problem}, also called the
{\em integer Carath\'eodory problem}, is the following: given a set $S\subseteq\Z^d$ and
a vector $a\in\Z^d$, decide if $a\in\mon(S)$, and if yes, find a decomposition
$a=\sum\l_k x^k$ with all $\l_k\in\Z_+$ and all $x^k\in S$.

Of course, the complexity of the problem depends on the presentation of $S$.
When $S$ is given explicitly as a set of vectors, this is simply integer programming,
so even the decision problem is already NP-complete. Here we consider the much more
difficult situation with $S$ given implicitly as the set $S=\{x\in\Z^d: Ax\leq b\}$
of integer points satisfying a given system of inequalities, where $A$ is a $c\times d$
integer matrix and $b\in \Z^c$. Under such a presentation, even for fixed dimension $d$,
the number of points in $S$ can be infinite or exponential in the encoding length of $A$ and $b$,
so it is unclear how to even write down in polynomial time an expression $a=\sum\l_k x^k$,
let alone find one. In spite of this, Eisenbrand and Shmonin showed in \cite{ES} that if
$a\in\mon(S)$ then there is an expression $a=\sum_{k=1}^m\l_k x^k$ with $m\leq2^d$.
Recently, Goemans and Rothvo{\ss} showed in \cite{GR}, using heavy machinery, that for fixed $d$,
the problem is polynomial-time solvable, with degree which is exponential in $d$.

Here we show that when $A$ is totally unimodular, which holds in the context
of the three-way table problem, the problem can be solved in polynomial time
even when $d$ is variable. (Note that $A$ being totally unimodular does not imply that so
is the matrix with columns in $S$, just take $d=1$ and $S=\{x\in\Z:2\leq x\leq 2\}=\{2\}$.)
In fact, we prove in Section 2 a more general result on the monoid problem for
sets $S=P\cap\Z^d$ for polyhedra $P$ in an oracle setup, and deduce the following corollary.

\vskip.2cm\noindent{\bf Corollary \ref{uic}}\hskip.2cm
{\em The monoid decomposition problem for any $S=\{x\in\Z^d: Ax\leq b\}$
and any $a$, with $A$ totally unimodular and $b$ integer, is solvable in time polynomial
in the binary-encoding length of $A,b$, and $a$, even when the dimension $d$ is variable.}

\vskip.3cm
We proceed as follows. In Section 2 we prove Corollary \ref{uic}. In Section 3 we discuss
multiway tables and multicommodity flows and use Corollary \ref{uic} to prove Theorem \ref{tables}
and Corollary \ref{commodity}. We conclude in Section 4 with some open problems.

\section{Unimodular integer Carath\'eodory}

As mentioned in the introduction, we solve here the monoid problem for a broad class of
sets of the form $S=P\cap\Z^d$ where $P$ are polyhedra presented by suitable oracles.
Corollary \ref{uic} will then follow as a special case. Throughout, all polyhedra $P\subset\R^d$
and all vectors $x\in\R^d$ are {\em rational}, and we will not indicate this further for brevity.
We will use the algorithmic theory of polyhedra developed in \cite{GLS}.
The {\em description complexity} of a polyhedron $P$ is the smallest
positive integer $\Delta$ such that $P$ admits a description $P=\{x\in\R^d: Ax\leq b\}$
with $A$, $b$ integer and $|A_{i,j}|,|b_i|\leq \Delta$ for all $i,j$. (We do not need to
know this description explicitly, and the number of inequalities may be exponential.)
A {\em separation oracle} for a polyhedron $P\subset\R^d$ is one that, queried on $x\in\R^d$,
either asserts that $x\in P$ or returns an $h\in\R^d$ such that $hy<hx$ for all $y\in P$.
In all algorithmic statements on polyhedra $P\subset\R^d$ involving oracles, an algorithm is
said to {\em run in polynomial time} if its running time including queries to the oracles
involved is polynomial in $d$, the binary-encoding $\log\Delta$ of the description complexity
of $P$, and other relevant inputs. See \cite{GLS} for more details.

\vskip.2cm
We begin with two simple lemmas.

\bl{scale}
Given polyhedron $P\subset\R^d$ presented by separation oracle and $a\in\R^d$, we can
in polynomial time either find $n\in\N$ with ${1\over n}a\in P$ or asserts none exists.
\el
\bpr
Using the separation oracle of $P$ it is possible to efficiently realize a separation oracle for
the intersection $Q:=P\cap\{\l a:0\leq\l\}$ of $P$ and the ray generated by $a$. Minimizing and
maximizing the linear function $ax$ over $Q$ using the algorithmic equivalence of separation
and optimization from \cite{GLS} we conclude with one of the following:
$Q=\emptyset$ so there is no $n$;
$Q=\{\l a:\alpha\leq\l\}$ and then if $\alpha>1$ then there is no $n$
whereas if $\alpha\leq 1$ then we can take $n=1$;
$Q=\{\l a:\alpha\leq\l\leq\beta\}$
and then if there is $n\in\N$ with ${1\over\beta}\leq n\leq{1\over\alpha}$ then we
can take it and otherwise there is no $n$.
\epr

A polyhedron $P\subset\R^d$ is {\em decomposable} if for every $n\in\N$ and every $x\in nP\cap\Z^d$
there are $x^1,\dots,x^n\in P\cap\Z^d$ with $x=x^1+\cdots+x^n$, where $nP:=\{ny:y\in P\}$.

Note that for any $S\subseteq\Z^d$, even $S=\emptyset$, we have that $a=0$ is trivially
in $\mon(S)$ with the empty decomposition. The next lemma deals with the case of $a\neq0$.

\bl{criterion}
Let $P\subset\R^d$ be a polyhedron and let $S:=P\cap\Z^d$. Let $a\in\Z^d$ be a nonzero vector.
If $a\in\mon(S)$ then there is an $n\in\N$ such that $a\in nP$. If $P$ is moreover
decomposable and there is an $n\in\N$ such that $a\in nP$ then $a\in\mon(S)$.
\el
\bpr
If $a\in\mon(S)$ then $a=\sum\l_kx^k$ with $\l_k\in\N$ and $x^k\in S\subseteq P$.
Let $n:=\sum\l_k\geq1$. Then $x:=\sum{\l_k\over n}x^k$ is a convex combination
of points in $P$ and therefore $x\in P$, so $a=nx\in nP$. If $P$ is decomposable
and $a\in nP$ for some $n\in\N$ then there are $x^1,\dots,x^n\in P\cap\Z^d=S$
with $a=x^1+\cdots+x^n$, so $a\in\mon(S)$.
\epr
A {\em decomposition oracle} for a decomposable $P$ is one that, queried on $n\in\N$ given in unary,
and on $x\in nP\cap\Z^d$, returns $x^1,\dots,x^n\in P\cap\Z^d$ with $x=x^1+\cdots+x^n$.

\vskip.2cm
We proceed to establish the efficient solution of the monoid problem over polyhedra defined
by oracles. For simplicity we provide the statement and proof for {\em pointed} polyhedra.
(A polyhedron is pointed if it has at least one vertex, which is equivalent to admitting an
inequality description with a matrix of full column rank.) The polyhedra appearing in typical
applications are indeed pointed. Moreover, in the specializations of the oracle result
to concrete polyhedra in Corollaries \ref{uic} and \ref{uic_extension} in the sequel,
we solve the monoid problem even for non pointed polyhedra.

\vskip.2cm
We will need the following result of \cite{GLS} (see Corollary 6.5.13 therein).
\bp{LP}
Given a pointed polyhedron $P\subset\R^d$ presented by a separation oracle and a point $x\in P$,
we can in polynomial time obtain vertices $x^0,\dots,x^k$ of $P$ for some $0\leq k\leq d$,
a point $y$ (possibly zero) in the recession cone of $P$, and positive rational numbers
$\l_0,\dots,\l_k$ satisfying $\sum_{i=0}^k\l_i=1$ and $x=y+\sum_{i=0}^k\l_ix^i$.
\ep

\vskip.2cm
We can now establish our oracle result.

\bt{uic_oracle}
The monoid decomposition problem over any set which is of the form $S:=P\cap\Z^d$ with
$P\subset\R^d$ any decomposable pointed polyhedron presented by a separation oracle and
endowed with a decomposition oracle is polynomial-time solvable.
\et
\bpr
Given any nonzero $a\in\Z^d$, we need to decide if $a\in\mon(S)$ and find a decomposition if yes.
We apply Lemma \ref{scale}. If there is no $n\in\N$ with $a\in nP$ then
$a\notin\mon(S)$ by Lemma \ref{criterion}. So assume we find $n\in\N$ with
$a\in nP$. Then $a\in\mon(S)$  by Lemma \ref{criterion} again.
We need to find a monoid decomposition of $a$.

But we {\em cannot} simply query the decomposition oracle on $n$ in unary and on $a$
to get the decomposition: the {\em crucial difficulty} is that the $n$ we got may be
very large, and only the {\em binary-encoding length} of $n$, not $n$ itself, is guaranteed
to be polynomial in the binary-encoding length of $a$ and the description complexity of $P$.

So instead we proceed as follows. We use Proposition \ref{LP} and obtain vertices $x^0,\dots,x^k$
of $P$ with $0\leq k\leq d$, point $y$ (possibly zero) in the recession cone of $P$, and positive
rational numbers $\l_0,\dots,\l_k$ with $\sum_{i=0}^k\l_i=1$ and ${1\over n}a=y+\sum_{i=0}^k\l_ix^i$.

Now, we claim that since $P$ is decomposable, its vertices are integer. Indeed, consider any
vertex $v$ of $P$. Since $P$ is rational so is $v$ and so for some $q\in\N$ we have
$qv\in qP\cap\Z^d$. Then $qv=z^1+\cdots+z^q$ for some $z^i\in P\cap\Z^d$ so
$v={1\over q}(z^1+\cdots+z^q)$ which implies $v=z^1=\cdots=z^q$ since $v$ is a vertex,
and therefore $v\in\Z^d$.

So the vertices $x^0,\dots,x^k$ that we obtained are in $P\cap\Z^d=S$. Now define
$${\bar n}\ :=\ \sum_{i=0}^k(n\l_i-\lfloor n\l_i\rfloor)
\ =\ n-\sum_{i=0}^k\lfloor n\l_i\rfloor\ ,$$
$${\bar a}\ :=\ ny+\sum_{i=0}^k(n\l_i-\lfloor n\l_i\rfloor)x^i
\ =\ a-\sum_{i=0}^k\lfloor n\l_i\rfloor x^i\ .$$
Suppose first that ${\bar n}=0$. Then $n\l_i$ is an integer for $i=0,\dots,k$ and therefore
$ny=a-\sum_{i=0}^k n\l_ix^i$ is an integer vector in the recession cone of $P$. Therefore
$x^0+ny\in P\cap\Z^d=S$. Now $n\l_0\geq 1$ and therefore we obtain the decomposition
$$a\ =\ (n\l_0-1)x^0+(x^0+ny)+\sum_{i=1}^k n\l_ix^i\ .$$
Next suppose ${\bar n}\neq0$. Then ${\bar n}$ is an integer satisfying $1\leq{\bar n}\leq d$.
Moreover, we have
$${1\over{\bar n}}{\bar a}\ =\ {n\over{\bar n}}y +
\sum_{i=0}^k {n\l_i-\lfloor n\l_i\rfloor\over{\bar n}}x^i\ .$$
So ${1\over{\bar n}}{\bar a}$ is the sum of a convex combination of vertices of $P$
and a vector in the recession cone of $P$, and hence is in $P$.
Therefore ${\bar a}\in{\bar n}P\cap\Z^d$. We now query the decomposition oracle of $P$ on ${\bar n}$
and ${\bar a}$ and obtain ${\bar a}=\sum_{i=1}^{\bar n}z^i$ for suitable $z^i\in P\cap\Z^d=S$.
This gives again a decomposition of $a$, and completes the proof,
$$a\ =\ \sum_{i=0}^k \lfloor n\l_i\rfloor x^i + \sum_{i=1}^{\bar n}z^i\ .\mepr$$

We next consider polyhedra defined by totally unimodular matrices. We need
an algorithmic version of the decomposition theorem of Baum and Trotter \cite{BT}.

\bl{decomposition}
For any totally unimodular matrix $A$ and any integer vector $b$, the polyhedron
$P:=\{x\in\R^d:Ax\leq b\}$ is decomposable. Moreover, there is a polynomial time
algorithm that, given such $A$ and $b$, realizes a decomposition oracle for $P$.
\el
\bpr
We show by induction on $n$ that given $n\in\N$ and $a\in nP\cap\Z^d$ we can find in
polynomial time $x^1,\dots,x^n\in P\cap\Z^d$ with $a=x^1+\cdots+x^n$. For $n=1$ simply
take $x^1=a$. Next consider $n>1$ and consider the following system in variable vector $x$,
\begin{equation*}\label{system0}
Ax \leq b\,,\ \ A(a-x) \leq (n-1)b\ .
\end{equation*}
Then $x:={1\over n}a$ is a real solution of this system, since
$a\in nP$ implies $A{1\over n}a\leq b$ and $A(a-{1\over n}a)=(n-1)A{1\over n}a\leq(n-1)b$.
Now, the defining matrix of the system \eqref{system0} consists of one block of $A$ and
one block of $-A$ and hence is totally unimodular since $A$ is, and the right hand side
of this system is integer. So the system also admits an integer solution $x^n\in\Z^d$
which can be found in polynomial time by linear programming.
Then $Ax^n\leq b$ and hence $x^n\in P\cap\Z^d$. Moreover, $A(a-x^n)\leq(n-1)b$
so $a-x^n\in(n-1)P\cap\Z^d$ and hence, by induction, we can find a decomposition
$a-x^n=\sum_{i=1}^{n-1}x^i$ with all $x^i\in P\cap \Z^d$.
This yields the decomposition $a=\sum_{i=1}^nx^i$.
\epr

\vskip.3cm
We can now conclude the following corollary mentioned in the introduction.

\bc{uic}
The monoid decomposition problem for any $S=\{x\in\Z^d: Ax\leq b\}$
and any $a$, with $A$ totally unimodular and $b$ integer, is solvable in time polynomial
in the binary-encoding length of $A,b$, and $a$, even when the dimension $d$ is variable.
\ec
\bpr
Let $A$, $b$, and $a\neq0$ be given input to the problem, so $P:=\{x\in\R^d:Ax\leq b\}$,
$S:=P\cap\Z^d$, and we need to decide if $a\in\mon(S)$ and find a decomposition if yes.
Since $P$ is not necessarily pointed, we proceed as follows. We let $T:=Q\cap\Z^d$ with
$$Q\ :=\ \{x\in\R^d\ :\ Ax\leq b\,,\ \ 0\leq\sign(a_i) x_i\leq|a_i|\,,\ \ i=1,\dots,d\}\ ,$$
where the {\em sign} of $r\in\R$ is $\sign(r):=1$ if $r\geq 0$ and $\sign(r):=-1$ if $r<0$.
Clearly $Q$ has a separation oracle and description complexity polynomial in the input.
Moreover, the system defining $Q$ is totally unimodular and hence
$Q$ is decomposable with a decomposition oracle by Lemma \ref{decomposition}.
We apply Lemma \ref{scale} to $P$. If there is no $n\in\N$ with $a\in nP$ then
$a\notin\mon(S)$ by Lemma \ref{criterion}. So assume we find $n\in\N$ with $a\in nP$.
Then also $a\in nQ\cap\Z^d$ and hence $a\in\mon(T)$ by Lemma \ref{criterion}.
Now $Q$ is a polytope hence pointed. So we can apply Theorem \ref{uic_oracle} to
$T=Q\cap\Z^d$ and obtain a decomposition $a=\sum\l_k x^k$ with $\l_k\in\Z_+$
and $x^k\in T\subseteq S$ as desired.
\epr

We conclude this section with an extension of Corollary \ref{uic} to the following
broader class of monoids. A {\em totally unimodular projection} is a polyhedron
of the form $P=\{x=Ly\in\R^d:y\in Q\}$ which is the linear projection of a polyhedron
$Q=\{y\in\R^c:Ay\leq b\}$, with $b$ an integer vector and $[A^{^T}\ L^{^T}]$ totally unimodular.
The special case $L=[I_d\ 0_{d\times(c-d)}]$ with $I_d$ the identity gives variable-eliminating
projections and the case $L=I_d$ and $c=d$ gives the polyhedra in Corollary \ref{uic}.

\vskip.2cm
We have the following extension of Corollary \ref{uic} to such polyhedra.
\bc{uic_extension}
The monoid decomposition problem over any totally unimodular\break projection $P$
can be solved in polynomial time even when the dimension $d$ is variable.
\ec
\bpr
Let $L$, $A$, $b$, and $a\neq 0$ be given input and let $S:=P\cap\Z^d$.
We need to decide if $a\in\mon(S)$ and find a decomposition if yes.
Note that the data gives separation oracles for both $P$ and $Q$
with description complexities polynomial in the input.

We apply Lemma \ref{scale} to $P$. If there is no $n\in\N$ with $a\in nP$ then $a\notin\mon(S)$ by
Lemma \ref{criterion}. So assume we find $n\in\N$ with $a\in nP$. Consider the system
\begin{equation*}\label{1}
Ay\ \leq\ nb\,,\quad Ly\ =\ a\ .
\end{equation*}
Since ${1\over n}a\in P$ there is a $y\in\R^c$ with $Ay\leq b$ and ${1\over n}a=Ly$.
Then $ny$ satisfies the system. Since $[A^{^T}\ L^{^T}]$ is totally unimodular,
we can find an integer solution $z$ to the system. So $z\in nQ\cap \Z^c$.
Let $T:=Q\cap\Z^c$. By Lemma \ref{criterion} we have that $z\in\mon(T)$.
Since $A$ is totally unimodular we can use Corollary \ref{uic} and find in polynomial time
a decomposition $z=\sum\l_kz^k$ for some $\l_k\in\Z_+$ and $z^k\in Q\cap\Z^c=T$. Let $x^k:=Lz^k$
for all $k$. Then $x^k\in P\cap\Z^d=S$ since $L$ is integer, and we obtain the decomposition
$$a\ =\ Lz\ =\ L\sum\l_kz^k\ =\ \sum\l_kLz^k\ =\ \sum\l_k x^k\ .\mepr$$

\section{Huge tables are fixed-parameter tractable}

As noted in the introduction, the three-way table problem is to decide if the following
set of nonnegative integer $l\times  m\times n$ tables is nonempty, and find a table if one exists,
$$\left\{x\in\Z_+^{l\times  m\times n}\ :\ \sum_i x_{i,j,k}=u_{j,k}
\,,\ \sum_j x_{i,j,k}=v_{i,k}\,,\ \sum_k x_{i,j,k}=w_{i,j}\right\}\ ,$$
with the line sums binary-encoded integers $u_{j,k}$, $v_{i,k}$, $w_{i,j}$
for $1\leq i\leq l$, $1\leq j\leq m$, and $1\leq k\leq n$, of binary-encoding length $\size(u,v,w)$.
This problem was shown in \cite{HOR} to be
fixed-parameter tractable when $l$ and $m$ are parameters, solvable in
time $O(f(l,m)\cdot n^3\cdot\size(u,v,w))$ for suitable computable
function $f(l,m)=(lm)^{O(lm)}$.

Regard now each table as a tuple $x=(x^1,\dots,x^n)$ consisting of
$n$ many $l\times m$ {\em layers}.
Following \cite{OS}, call the problem {\em huge} if the
variable number $n$ of layers is encoded in {\em binary}.
We are then given $t$ {\em types} of layers, where each type $k$ has its column sums
vector $u^k\in\Z_+^m$ and row sums vector $v^k\in\Z_+^l$. In addition, we are given positive
integers $n_1,\dots,n_t,n$ with $n_1+\cdots+n_t=n$, all encoded in binary. A feasible
table $x=(x^1,\dots,x^n)$ then must have first $n_1$ layers of type $1$, next $n_2$ layers
of type $2$, and so on, with last $n_t$ layers of type $t$. The special case of $t=1$ type
is the {\em symmetric case}, where all layers have the same row and column sums.

\bt{tables} The huge $l\times m\times n$ table problem with $t$
types, parameter $l$, and $n$ variable and binary-encoded, is
fixed-parameter tractable in the following situations:
\begin{enumerate}
\item
when $m$ is also a parameter and $t$ is variable and unary-encoded;
\item
when $t$ is also a parameter and $m$ is variable and unary-encoded.
\end{enumerate}
\et

\vskip.2cm
\bpr
We first formulate the symmetric case as a monoid decomposition problem.

Let
$$S\ :=\ \left\{z\in\Z_+^{lm}\cong\Z_+^{l\times m}\ :\ \sum_i z_{i,j}=u_j
\,,\ \sum_j z_{i,j}=v_i\right\}\ =\ \{z\in\Z_+^{d}\, :\, Az=b\}$$
with $d=lm$, suitable $b$, and $A$ the $(l+m)\times lm$ vertex-edge incidence matrix of the
complete bipartite graph $K_{l,m}$ which is well known to be totally unimodular.

Note that even when $l$ and $m$ are fixed, the number of elements of $S$ is typically
exponential in the binary-encoding length $\size(u,v)$ of the row and column sums.

Now, if $x=(x^1,\dots,x^n)$ is a feasible symmetric table then $x^k\in S$ for all $k$ and
$\sum_{k=1}^n x^k=w$ is the vertical line sum vector so $Aw=\sum_{k=1}^n Ax^k =nb$.
So assume $Aw=nb$, which is easy to check, else there is no feasible table and we are done.
Assume also $b\neq0$ else the unique feasible table is zero and we are done again.
Then $w\in\mon(S)$ by Lemmas \ref{criterion} and \ref{decomposition},
so $w=\sum\l_k z^k$ for some $\l_k\in\Z_+$ and $z^k\in S$.

We then have
$$\sum\l_k b\ =\ \sum\l_k Az^k\ =\ A\sum\l_k z^k\ =\ Aw\ =\ nb\ ,$$
so $\sum\l_k=n$. So there is a feasible table $x=(x^1,\dots,x^n)$
with $\l_k$ layers equal to $z^k$ for all $k$. By Corollary \ref{uic} we can
solve this monoid decomposition problem in time polynomial in $l$ and $m$ and $\size(u,v,w)$
and find the $\l_k\in\Z_+$ and $z^k\in S$ which provide a compact representation
of a feasible huge symmetric three-way table $x$.

We proceed to the general case of huge tables with $t$ types.
The solution has two steps. First, a regular (non huge) compressed problem
over $l\times m\times t$ tables is derived from the huge problem data, with the same
vertical sums $w$, and column sums $n_ku^k$ and row sums $n_kv^k$ for $k=1,\dots,t$.
We solve the compressed problem by the Graver bases
methods of \cite{HOR}, either in time $O(f(l,m)\cdot t^3\cdot\size(n_ku^k,n_kv^k,w))$,
or in time $O(f(l,t)\cdot m^3\cdot\size(n_ku^k,n_kv^k,w))$, according to which
of $m$ and $t$ is chosen to be the parameter and which is chosen to be variable.
If the original problem has a feasible table $x=(x^1,\dots,x^n)$ then,
taking $y^1$ to be the sum of the first $n_1$ layers of $x$, taking $y^2$ to
be the sum of the next $n_2$ layers of $x$, and so on, with lastly taking $y^t$
to be the sum of the last $n_t$ layers of $x$, we obtain a table $y=(y^1,\dots,y^t)$ which is
feasible in the compressed problem. So assume we found a table $y$ which is feasible
in the compressed problem, else the original problem is infeasible and we are done.

Second, for $k=1,\dots,t$ we consider the huge symmetric table problem which asks for
a huge $l\times m\times n_k$ table with vertical sums given by $y^k$ and with column
sums $u^k$ and row sums $v^k$. Each of these $t$ huge symmetric table problems
is now formulated as a monoid problem as just explained above with the same matrix $A$
and a suitable $b^k$ defined from $u^k$ and $v^k$. Since $y^k$ has
column sums $n_ku^k$ and row sums $n_kv^k$, we have $Ay^k=n_kb^k$ and so,
as explained above, this symmetric problem is feasible and by Corollary \ref{uic}
we can find in polynomial time a compact representation of a feasible
$l\times m\times n_k$ table. The concatenation of these compact representations
provides a compact representation of an $l\times m\times (n_1+\cdots+n_t)=l\times m\times n$
table which provides the desired solution of the original huge table problem.
\epr

This theorem also has a consequence to the following huge
multicommodity flow problem over the complete bipartite graph.
There are $l$ commodities, $m$ suppliers, and nonnegative integer numbers $s^k_i$ of
units that supplier $i$ is to supply of commodity $k$. There are $t$ consumer types, where,
for $r=1,\dots,t$, we have $n_r$ consumers of type $r$, nonnegative integer numbers $c^k_r$
of units that each consumer of type $r$ is to consume of commodity $k$, and nonnegative
integer capacities $u_{i,r}$ of allowed flow of all commodities from supplier $i$ to each
consumer of type $r$. Adding one slack commodity, we may assume that the capacities should be
attained with equality. The problem is to find a (compact representation of a) feasible flow
$x^k_{i,j}$ from each supplier $i$ to each consumer $j$ of each commodity $k$. The numbers
$n_k$ of consumers of each type are encoded in binary, so there is a huge number $n_1+\cdots+n_t$
of consumers. It is then not hard to see that this can be directly encoded as a suitable
huge three-way table problem. Theorem \ref{tables} then implies the following statement.

\bc{commodity}
The huge multicommodity flow problem with $l$ commodities, $t$
consumer types, unary-encoded number $m$ of suppliers, binary-encoded
supplies, consumptions and capacities $s^k_i$, $c^k_r$ and $u_{i,r}$,
and binary-encoded numbers $n_r$ of consumers of type $r=1,\dots,t$,
parameterized by $l$ and $t$, is fixed-parameter tractable.
\ec

We conclude this section with an extension of Theorem \ref{tables} to
a class of {\em huge $n$-fold integer programming problems} defined as follows.
The {\em $n$-fold product} of a $c\times d$ integer matrix $A$
is the following $(d+cn)\times(dn)$ matrix, with $I_d$ the identity,
$$A^{[n]}\quad:=\quad
\left(
\begin{array}{cccc}
  I_d    & I_d    & \cdots & I_d    \\
  A    & 0      & \cdots & 0      \\
  0      & A    & \cdots & 0      \\
  \vdots & \vdots & \ddots & \vdots \\
  0      & 0      & \cdots & A    \\
\end{array}
\right)\quad .
$$
The {\em $n$-fold integer programming feasibility problem} is to decide if the following set,
\begin{equation}\label{classical-n-fold}
\left\{x\in\Z^{dn}\ :\ A^{[n]}x=b\,,\ l\leq x\leq u\right\}\ ,
\end{equation}
is nonempty, and find a feasible point if there is one, where $b\in\Z^{d+cn}$
and $l,u\in\Z^{dn}$. See \cite{Onn} for more details and for the many
applications of this class of problems.

In \cite{HOR} it was shown to be fixed-parameter tractable parameterized by $d$
and $a\in\N$ an upper bound on $|A_{i,j}|$ for all $i,j$, solvable in time
$O(f(a,d)\cdot n^3\cdot\size(b,l,u))$, with $\size(b,l,u)$ the bit size of $b,l,u$,
and $f(a,d)=(ad)^{O(d^3)}$ a computable function.

The vector ingredients of an $n$-fold program are naturally arranged in {\em bricks},
with $x=(x^1,\dots,x^n)$ and likewise $l$ and $u$ with $x^i,l^i,u^i\in\Z^d$
for $i=1,\dots,n$, and with $b=(b^0,b^1,\dots,b^n)$ with $b^0\in\Z^d$
and $b^i\in\Z^c$ for $i=1,\dots,n$. Following \cite{OS}, the $n$-fold program
is called {\em huge} if $n$ is encoded in {\em binary}. More precisely,
we are now given $t$ {\em types} of bricks, where each type $k=1,\dots,t$ has
its lower and upper bounds $l^k,u^k\in\Z^d$ and right-hand side $b^k\in\Z^c$.
A brick $x^i\in\Z^d$ has type $k$ if $Ax^i=b^k$ and $l^k\leq x^i\leq u^k$.
Also given are $b^0\in\Z^d$ and $n_1,\dots,n_t,n\in\N$ with $n_1+\cdots+n_t=n$, all encoded in
binary. A feasible point $x=(x^1,\dots,x^n)$ now must have first $n_1$ bricks of type $1$,
next $n_2$ bricks of type $2$, and so on, with last $n_t$ bricks of type $t$, and also
satisfy $\sum_{i=1}^n x^i=b^0$. When the defining matrix $A$ is totally unimodular, which in
particular implies $a=1$ holds, a proof similar to that of Theorem \ref{tables}, the details of
which are omitted, where bricks replace layers, leads to the following theorem.

\bt{nfold}
The huge $n$-fold integer programming problem with $t$ types over any totally unimodular
$c\times d$ matrix $A$, parameterized by $d$, with $t$ unary-encoded and with $b^0$, $n$,
$b^k,l^k,u^k$, and $n_k$ binary-encoded for $k=1,\dots,t$, is fixed-parameter tractable.
\et

\section{Open problems}

We now raise several remaining open problems. First, the complexity
of huge tables of higher dimensions is unsettled. For such tables,
the layers are tables of dimension at least three, and therefore the
matrix which defines the resulting monoid is no longer totally
unimodular. In particular, what is the complexity of deciding the
existence of a huge four-way $k\times l\times m\times n$ table with
given sums, with $k,l,m$ fixed and $n$ encoded in binary, with $t$
types? It is known that for fixed $t$ the problem is in P, and for
variable $t$ it is in NP intersect coNP but is not known to be in P
even for $3\times 3\times 3\times n$ tables, see \cite{OS}. We also
do not know whether the problem, with $k,l,m$ as parameters, with
variable or even fixed $t$, is fixed-parameter tractable.

Next, we discuss bin packing. We need to pack items of $d$ types in identical bins.
Each item of type $i$ has a positive integer volume $v_i$ and there are $n_i$
items of type $i$ to be packed. Each bin has a positive integer volume $v$.
The question is, given $n$, whether $n$ bins suffice to pack all items.
(The minimal possible number of bins needed can then be found by binary search.)
We consider the huge version of the problem, usually referred to as the
{\em cutting stock problem}, where all data, including the numbers $n_i$ of items
of each type $i$, and $n$, are encoded in binary. The study of this huge version
goes back to the classical paper \cite{GG} by Gilmore and Gomory.
We formulate this problem as a monoid decomposition problem in $\Z^{d+1}$ as follows. Let
$$S\ :=\ \{z\in\Z_+^{d+1}\ :\ z_0=1\,,\ \sum_{i=1}^d v_iz_i\leq v\}
\ =\ \{z\in\Z_+^{d+1}\, :\, z_0=1\,,\ Az\leq b\}$$ with
$A=(0,v_1,\dots,v_d)$, $b=v$, and variables $z=(z_0,z_1,\dots,z_d)$.
Then $z\in S$ if and only if $z_0=1$ and $(z_1,\dots,z_d)$ is an
{\em admissible packing pattern} which means that it is possible to
pack $z_i$ items of type $i$ for $i=1,\dots,d$ in a single bin. Now
let $a:=(n,n_1,\dots,n_d)$. Then $a=\sum\l_k z^k$ with $\l_k\in\Z_+$
and $z^k\in S$ if and only if there is a packing of all items using
$\sum\l_k=n$ bins, where for each $k$ there are $\l_k$ bins packed
in pattern $(z^k_1,\dots,z^k_d)$. The solution of this monoid
decomposition problem allows to decide if there is a packing with
$n$ bins and if there is, to find it. McCormick, Smallwood and
Spieksma showed in \cite{MSS} that the problem is polynomial-time
solvable for $d=2$ types, and asked about higher $d$. This was
resolved only very recently by Goemans and Rothvo{\ss} in \cite{GR}
who showed that it can be solved in polynomial time for any fixed $d$,
as a consequence of their solution of the monoid problem for fixed $d$.
Unfortunately, the matrix $A$ above is not totally unimodular and therefore
Corollary \ref{uic} does not apply. So it remains open whether the cutting
stock problem, with the number $d$ of types as a parameter, is
fixed-parameter tractable or not.

Finally, it remains an important open question whether the monoid
problem for general matrices, parameterized by the dimension $d$, is
fixed-parameter tractable.

\section*{Acknowledgment}

This research was partially supported by the Dresner Chair at the Technion.
I am indebted to the referees for simplifying and strengthening the results in Section 2.


\begin{thebibliography}{}

\bibitem{BT}
Baum, S., Trotter, L.E., Jr.:
Integer rounding and polyhedral decomposition for totally unimodular systems.
Lecture Notes in Economical and Mathematical Systems 157:15--23 (1978)

\bibitem{BO}
Berstein, Y., Onn, S.:
The Graver complexity of integer programming.
Annals of Combinatorics 13:289--296 (2009)

\bibitem{DHOW}
De Loera, J., Hemmecke, R., Onn, S., Weismantel, R.:
$n$-Fold integer programming.
Discrete Optimization 5:231--241 (2008)

\bibitem{DO1}
De Loera, J., Onn, S.:
The complexity of three-way statistical tables.
SIAM Journal on Computing 33:819--836 (2004)

\bibitem{DO2}
De Loera, J., Onn, S.:
All linear and integer programs are slim 3-way transportation programs.
SIAM Journal on Optimization 17:806--821 (2006)

\bibitem{DF}
Downey, R.G., Fellows, M.R.:
Fundamentals of Parameterized Complexity. Texts in Computer Science, Springer (2013)

\bibitem{ES}
Eisenbrand, F., Shmonin, G.:
Carath\'eodory bounds for integer cones.
Operations Research Letters 34:564--568 (2006)

\bibitem{FR}
Fienberg, S.E., Rinaldo, A.:
Three centuries of categorical data analysis:
Log-linear models and maximum likelihood estimation.
Journal of Statistical Planning and Inference 137:3430--3445 (2007)

\bibitem{GG}
Gilmore, P.C., Gomory, R.E.:
A linear programming approach to the cuttingstock problem.
Operations Research 9:849--859 (1961)

\bibitem{GR}
Goemans, M.X., Rothvo\ss, T.:
Polynomiality for bin packing with a constant number of item types.
In: Proceedings of the Symposium on Discrete Algorithms 25:830--839 (2014)

\bibitem{GLS}
Gr\"otschel, M., Lov\'asz, L., Schrijver, A.:
Geometric Algorithms and Combinatorial Optimization.
Second edition (1993), Springer

\bibitem{HOR}
Hemmecke, R., Onn, S., Romanchuk, L.:
$n$-Fold integer programming in cubic time.
Mathematical Programming 137:325--341 (2013)

\bibitem{KT}
Kudo, T., Takemura, A.:
A lower bound for the Graver complexity of the incidence matrix of a complete bipartite graph.
Journal of Combinatorics 3:695--708 (2012)

\bibitem{MSS}
McCormick, S.T., Smallwood, S.R., Spieksma, F.C.R.:
A polynomial algorithm for multiprocessor scheduling with two job lengths.
Mathematics of Operations Research 26:31--49 (2001)

\bibitem{Mot}
Motzkin, T.S.:
The multi-index transportation problem.
Bulletin of the American Mathematical Society 58:494 (1952)

\bibitem{Onn}
Onn, S.: Nonlinear Discrete Optimization.
Zurich Lectures in Advanced Mathematics,
European Mathematical Society (2010),
available online at: {\tt http://ie.technion.ac.il/$\sim$onn/Book/NDO.pdf}

\bibitem{OS}
Onn, S., Sarrabezolles, P.: Huge unimodular $n$-fold programs.
SIAM Journal on Discrete Mathematics 29:2277--2283 (2015)

\end{thebibliography}
\end{document}